\documentclass{article}
\usepackage{amsfonts,amsmath,amsthm,graphicx,float,stix}
\usepackage[font=large,labelfont=bf]{caption}
\usepackage[hidelinks]{hyperref}
\usepackage{apacite}
\theoremstyle{definition}
\newtheorem{definition}{Definition}
\newtheorem{example}{Example}
\newtheorem{theorem}{Theorem}

\newcommand{\set}{\emph{SET}}
\newcommand{\stun}{\emph{STUN}}
\DeclareMathOperator{\Aut}{Aut}

\begin{document}
\title{Toward a Combinatorial Theory of \set\\ and Related Card Games}
\author{Jonathan Schneider \small(Colby College)}
\date{First draft: June 14, 2019\\ This draft: October 25, 2021}
\maketitle

\section{The \set\, hypercube}\label{intro}
The popular card game \set\footnote{Published by Set Enterprises, first published in 1991.} features a deck of 81 distinct cards, each one distinguished by its unique combination of color, shape, fill, and number of symbols.  Each card has one of three possible values in each of these four attributes.  The goal of the game is to select three cards (called a ``Set'') that, in each attribute, either all match or all differ.  \cite{set}

\begin{figure}[H]\centering
\includegraphics[width=0.5in]{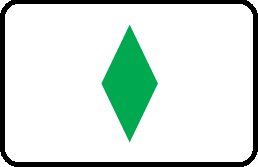}
\includegraphics[width=0.5in]{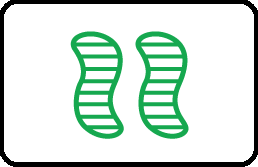}
\includegraphics[width=0.5in]{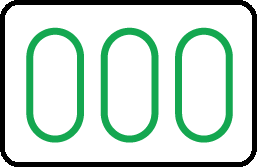}
\caption{A ``Set''.}\label{aset}
\end{figure}

\begin{figure}\centering
\includegraphics[height=.4\textheight]{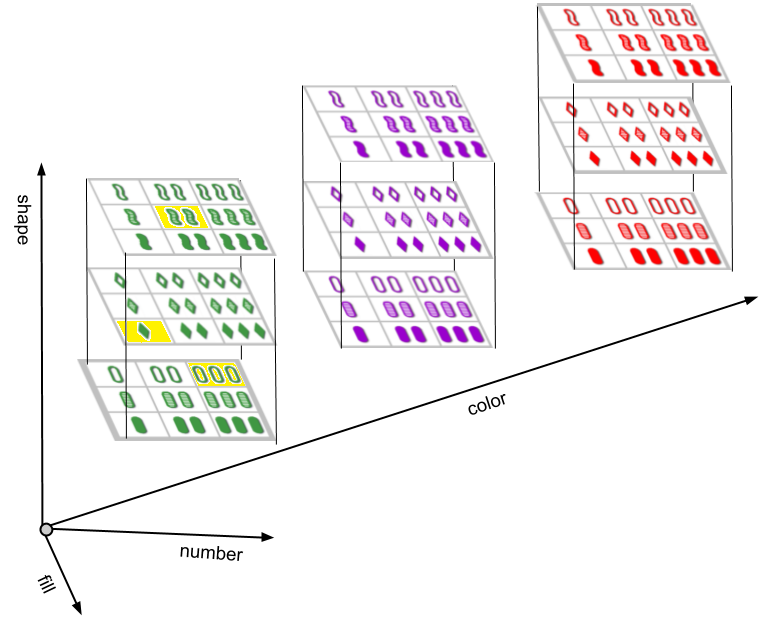}
\caption{The \set\, deck, with a ``Set'' highlighted in yellow.}\label{setincube}
\end{figure}

The deck can be plotted as a $3\times3\times3\times3$ hypercube lattice (figure \ref{setincube}).  A ``Set'' is distributed through this lattice so that no two of its cards occupy the same level along any axis, unless the third card is there too.  The ``Set'' of figure \ref{aset} is highlighted yellow in figure \ref{setincube}.  Notice that the three cards occupy the same level on the color axis but are spread across all three levels along the shape, fill, and number axes.

\begin{figure}\centering
\includegraphics[height=.4\textheight]{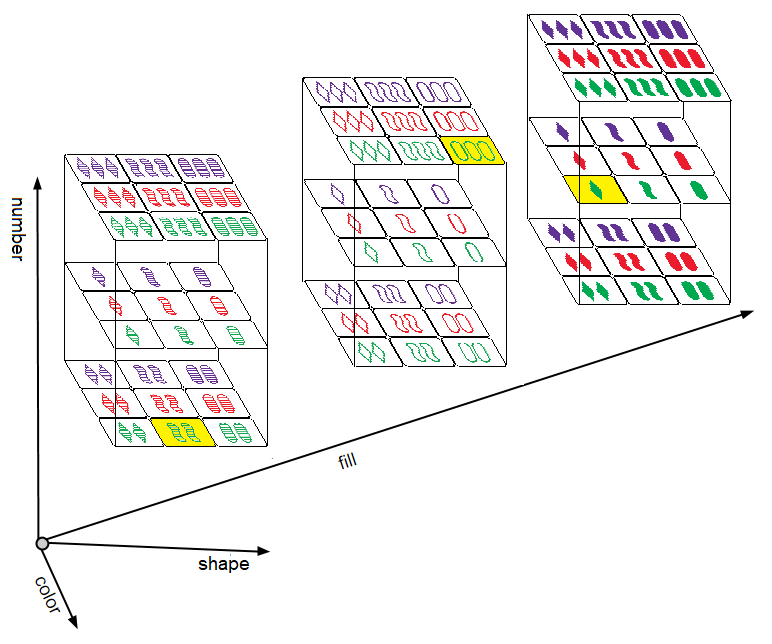}
\caption{The \set\, deck, with its attributes and values permuted.}\label{permutedcube}
\end{figure}

The order of the axes, and the order of values along each axis, were chosen arbitrarily for figure \ref{setincube}. The deck could have alternatively been drawn as in figure \ref{permutedcube}.  Whereas in figure \ref{setincube}, the axes are fill-number-color-shape (in that order), in figure \ref{permutedcube} the axes are color-shape-fill-number.  The values along the color axis in figure \ref{setincube} are green-purple-red (in that order), whereas in figure \ref{permutedcube} they are purple-red-green; the other attributes have their values reordered similarly.  The ``Set'' of figure \ref{aset} is still highlighted yellow.  Its cards are still aligned on one axis and spread out along the other three.  

A number of facts are immediately clear.  First of all, any two cards are part of a unique ``Set'', the third card being determined.  From this one can infer that the top three cards of a shuffled deck are a ``Set'' 1/79 of the time, and that there are ${81\choose3}/79=1080$ different ``Sets'' possible, out of ${81\choose3}=85,320$ card-triples that exist.

\begin{figure}\centering
Zero common attributes: $\vcenter{\hbox{
\includegraphics[width=0.4in]{16}
\includegraphics[width=0.4in]{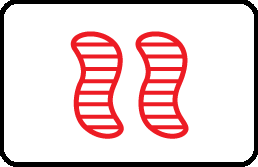}
\includegraphics[width=0.4in]{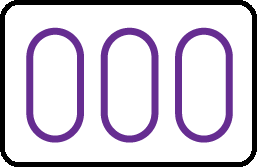}}}$\\
\vspace{4pt}
One common attribute: $\vcenter{\hbox{
\includegraphics[width=0.4in]{16}
\includegraphics[width=0.4in]{35}
\includegraphics[width=0.4in]{81}}}$\\
\vspace{4pt}
Two common attributes: $\vcenter{\hbox{
\includegraphics[width=0.4in]{16}
\includegraphics[width=0.4in]{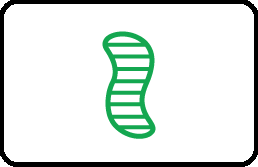}
\includegraphics[width=0.4in]{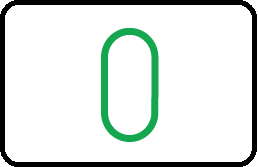}}}$\\
\vspace{4pt}
Three common attributes: $\vcenter{\hbox{
\includegraphics[width=0.4in]{16}
\includegraphics[width=0.4in]{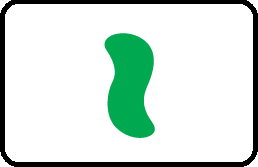}
\includegraphics[width=0.4in]{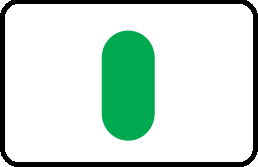}}}$\\
\caption{Four types of ``Set''.}\label{foursets}
\end{figure}

Secondly, ``Sets'' come in four types, depending on how many attributes the cards match in.  The ``Set'' of figure \ref{aset} has one attribute matching (color) and three attributes differing (shape, fill, \& number).  Examples of all four types of ``Set'' are shown in figure \ref{foursets}.

\begin{figure}\centering
\includegraphics[width=\textwidth]{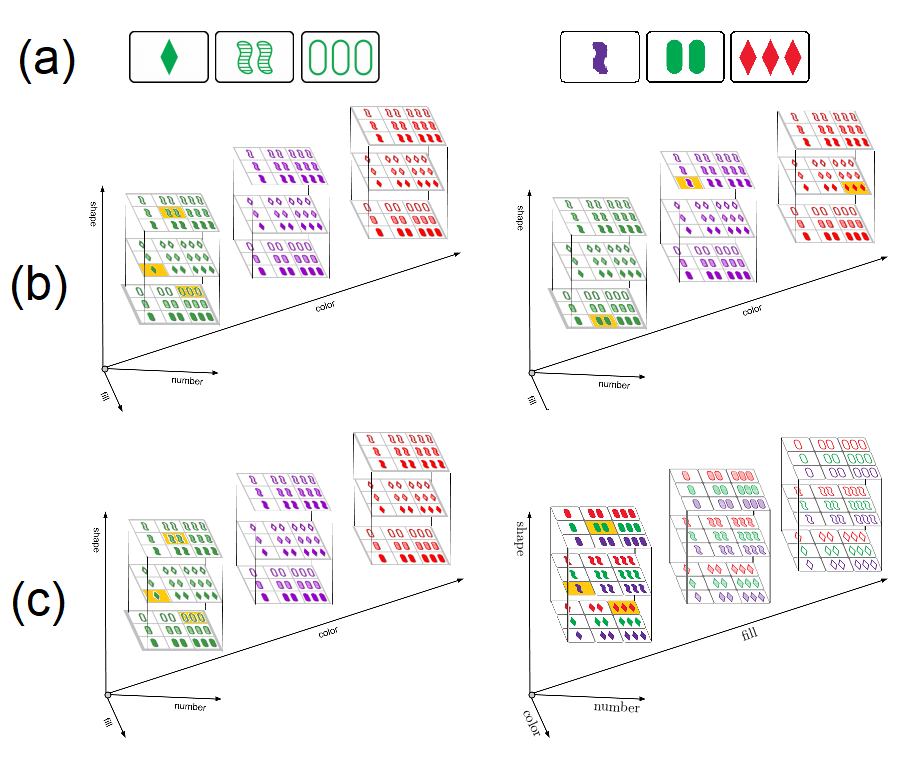}
\caption{\\(a) Two ``Sets'', both having one common attribute.\\ (b) The two ``Sets'' are highlighted in the same hypercube.\\ (c) The right-hand hypercube has been permuted so that the position of the highlighted ``Set'' matches the left-hand figure.}\label{similarsets}
\end{figure}

In figure \ref{similarsets}, we compare two ``Sets'', both having one common attribute.  Superficially, the two ``Sets'' are arranged in different positions within the hypercube.  However, by permuting the axes and values in the diagram, they can be put into correspondence.

\begin{figure}\centering
\includegraphics[width=\textwidth]{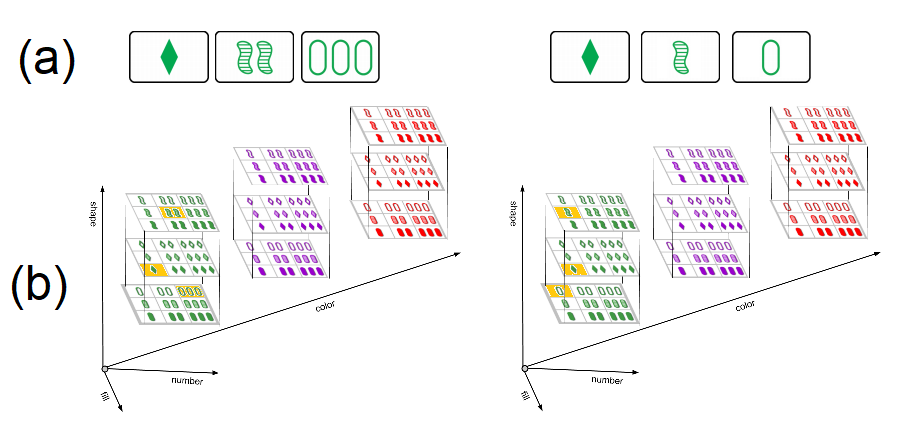}
\caption{\\(a) Two ``Sets'', having different numbers of common attributes.\\ (b) The two ``Sets'' are highlighted in the same hypercube.  No permutation of either cube can bring the yellow cards' positions into correspondence.}\label{distinctsets}
\end{figure}

In figure \ref{distinctsets}, we again compare two ``Sets'', but this time one ``Set'' has one common attribute while the other has two.  No permutation of the hypercube axes or values can put their arrangements into correspondence.  These ``Sets'' are \emph{fundamentally} different in this regard.

This suggests a general principle for comparing two collections of cards (or ``hands'', as we'll call them).  If two hands can be put into bijective correspondence by some permutation of the four attributes and the three values of each attribute, we regard them as hands \emph{of the same type}.  Among ``Sets'', in particular, there are four types (figure \ref{foursets}), but there are also many other types of hand with three or more cards (see, for example, figure \ref{tabthreehands}).

The purpose of this paper is to explore this equivalence relation on subsets of the \set\, deck, and to propose other games that can be invented on that basis.  In the process, we will enumerate some of the equivalence classes and their elements.

\pagebreak
\section{Isomorphism of hands}

The distinction between the four ``Sets'' of figure \ref{foursets} can be generalized to any collections of cards.

\begin{definition}  
A subset of the 81-card \set\, deck is called a \emph{hand}.  We fix the following symbols:  $A$ is the set of attributes, $A=\{$color, shape, fill, number$\}$. $V_a$ is the set of values for attribute $a\in A$.  Specifically,
\begin{align*}
V_\text{color}&=\{\text{red, green, purple}\},\\
V_\text{shape}&=\{\text{oval, diamond, squiggle}\},\\
V_\text{fill}&=\{\text{solid, empty, stripe}\},\\
V_\text{number}&=\{\text{single, double, triple}\}.
\end{align*}
If $x$ is a card, $v_a(x)$ is the card's value for attribute $a$.  For example,
\[
v_\text{color}\left(\vcenter{\hbox{\includegraphics[height=0.2 in]{16}}}\right)=\text{green.}
\]
\end{definition}
\begin{minipage}{\textwidth}\begin{definition}\label{defiso}
Two hands $H$ and $H'$ are \emph{isomorphic} if there exist
\begin{itemize}
	\item a bijection between the hands, $\varphi:H\to H'$
	\item a permutation of the four attributes, $\psi:A\to A$
	\item for each attribute $a\in A$, a bijection $\vartheta_a: V_a\to V_{\psi a}$
\end{itemize}
such that, for each attribute $a\in A$, the following diagram commutes:
\[
\includegraphics[width=1.5in]{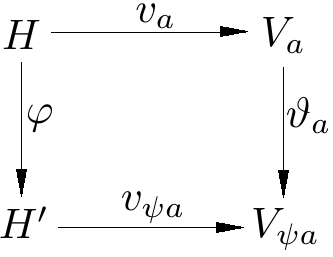}
\]
In this case, $\varphi$ is the \emph{isomorphism} induced by $\psi$ and the four $\vartheta_a$'s.
\end{definition}\end{minipage}
\begin{example}\label{exiso}
Let $H$ and $H'$ be the four-card hands shown.  We claim that $H$ and $H'$ are isomorphic.
\[
\begin{matrix}
H=	&\vcenter{\hbox{\includegraphics[width=.4in]{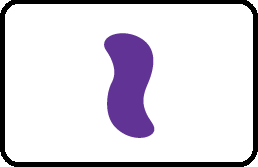}}}
	&\vcenter{\hbox{\includegraphics[width=.4in]{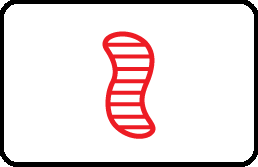}}}
	&\vcenter{\hbox{\includegraphics[width=.4in]{81}}}
	&\vcenter{\hbox{\includegraphics[width=.4in]{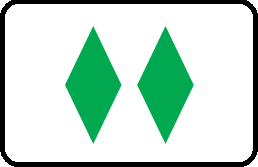}}}\\
\varphi:&\mapsdown&\mapsdown&\mapsdown&\mapsdown\\
H'=	&\vcenter{\hbox{\includegraphics[width=.4in]{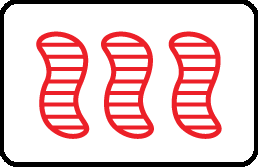}}}
	&\vcenter{\hbox{\includegraphics[width=.4in]{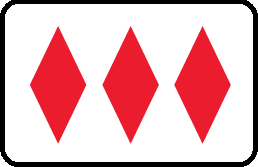}}}
	&\vcenter{\hbox{\includegraphics[width=.4in]{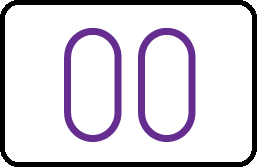}}}
	&\vcenter{\hbox{\includegraphics[width=.4in]{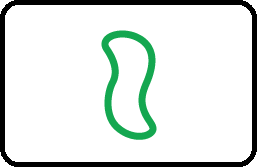}}}
\end{matrix}
\]
The indicated bijection $\varphi:H\to H'$ is induced by the following permutation of attributes and values.  (Verify that each card in $H$ is transformed accordingly.)
\begin{center}\begin{tabular}{|rcl|rcl|}
	\hline
	\textbf{COLOR}	&$\Longmapsto$& \textbf{FILL}	&\textbf{FILL}	&$\Longmapsto$& \textbf{SHAPE}\\
	red				&$\mapsto$& 	solid 			& empty			&$\mapsto$& 	oval\\
	green			&$\mapsto$& 	empty			& solid			&$\mapsto$& 	squiggle\\
	purple			&$\mapsto$& 	stripe			& stripe			&$\mapsto$& 	diamond\\
	\hline
	\textbf{NUMBER}	&$\Longmapsto$&\textbf{NUMBER}&\textbf{SHAPE}	&$\Longmapsto$&\textbf{COLOR}\\
	single			&$\mapsto$& 	triple 			& diamond		&$\mapsto$& 	green\\
	double			&$\mapsto$& 	single			& oval			&$\mapsto$& 	purple\\
	triple			&$\mapsto$& 	double			& squiggle		&$\mapsto$& 	red\\
	\hline	
\end{tabular}\end{center}
\noindent The large arrows ($\Longmapsto$) in this table are the permutation $\psi:A\to A$, and the small arrows ($\mapsto$) are the bijections $\vartheta_a:V_a\to V_{\psi a}$, as in definition 2. The upper-left box in the table gives rise to the commuting square in figure \ref{colorsquare}; the other three boxes also give commuting squares (exercise).
\begin{figure}\centering
\includegraphics[width=.7 \textwidth]{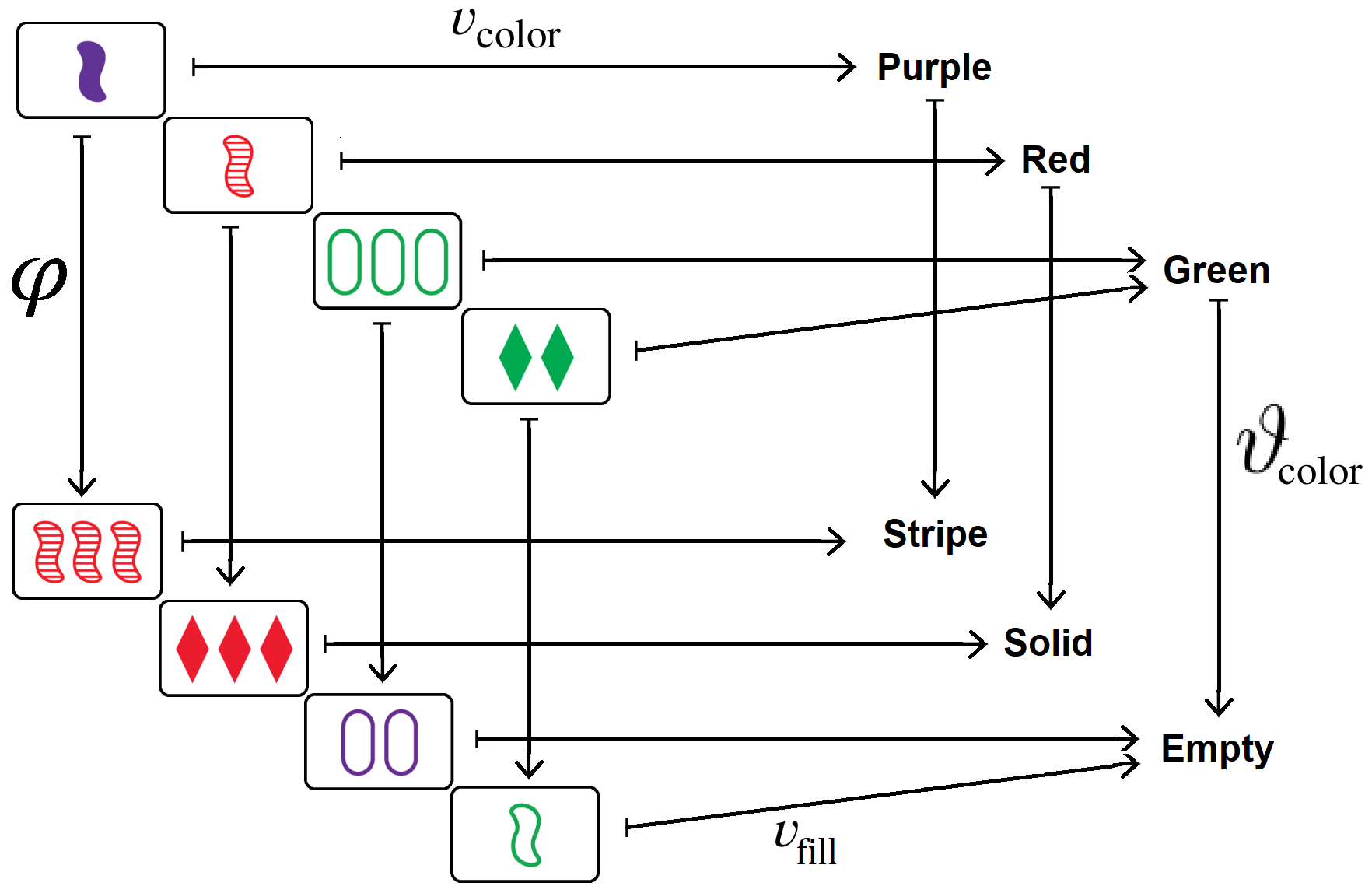}
\caption{One of four commuting squares used in example \ref{exiso}.}\label{colorsquare}
\end{figure}

We have determined that the hands $H$ and $H'$ are isomorphic by producing bijections $\varphi,\psi,\vartheta_a$ satisfying definition \ref{defiso}.  In general, searching for the bijections to induce a desired isomorphism is a lengthy (but finite) process.  Luckily, it is not necessary to determine the value maps $\vartheta_a$ completely.  Once $\varphi$ and $\psi$ have been declared, we can simply check that the hands \emph{split} in corresponding ways for each attribute.  In this example, $H$ splits according to color and $H'$ splits according to fill like this:
\[
\begin{matrix}
H=&\bigg(\vcenter{\hbox{\includegraphics[width=.4in]{4}}}\bigg)
	&\bigg(\vcenter{\hbox{\includegraphics[width=.4in]{28}}}\bigg)
	&\bigg(\vcenter{\hbox{\includegraphics[width=.4in]{81}}}
	&\vcenter{\hbox{\includegraphics[width=.4in]{17}}}\bigg)\\
\varphi:&\mapsdown&\mapsdown&\mapsdown&\mapsdown\\
H'=&\bigg(\vcenter{\hbox{\includegraphics[width=.4in]{30}}}\bigg)
	&\bigg(\vcenter{\hbox{\includegraphics[width=.4in]{12}}}\bigg)
	&\bigg(\vcenter{\hbox{\includegraphics[width=.4in]{77}}}
	&\vcenter{\hbox{\includegraphics[width=.4in]{61}}}\bigg)
\end{matrix}
\]
Since $\varphi$ respects this splitting for each attribute, it is an isomorphism.  Thus, the problem of deciding equivalence reduces to checking correspondence in how the hands split via each attribute.
\end{example}

\begin{example}
Let $H$ and $H'$ be the four-card hands shown.  Are they isomorphic?
\begin{align*}
	H=&\;\;\vcenter{\hbox{\includegraphics[width=.4in]{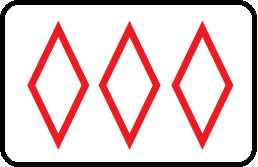}
		\includegraphics[width=.4in]{78}
		\includegraphics[width=.4in]{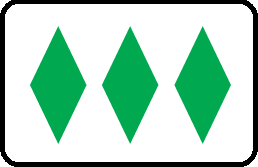}
		\includegraphics[width=.4in]{30}}}\\[.5em]
	H'=&\;\;\vcenter{\hbox{\includegraphics[width=.4in]{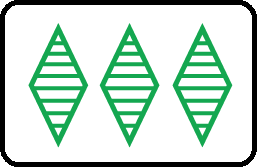}
		\includegraphics[width=.4in]{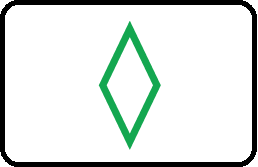}
		\includegraphics[width=.4in]{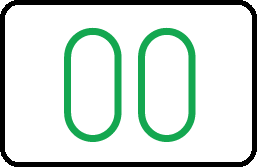}
		\includegraphics[width=.4in]{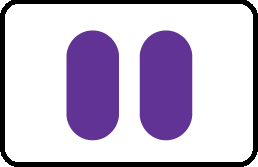}}}
\end{align*}
No, they are not.  In principle, one could decide this by checking every possible set of bijections $\psi$ and $\vartheta_a$ (only finitely many exist) to verify that none of them induces a bijection $\varphi:H\to H'$.    However, one can reach this conclusion much faster by observing how the hands split via each attribute.  For instance, $H'$ splits via color into 3 and 1 cards, whereas $H$ does not split into 3 and 1 via any attribute.  Alternatively, one could observe that $H'$ includes a pair of cards with no common attributes, whereas $H$ does not.
\end{example}

Isomorphism, as defined in this section, is the same as the hypercube correspondence described in section 1.  Permuting the hypercube axes is now $\psi:A\to A$; scrambling the order of the values along each axis is now $\vartheta_a:V_a\to V_{\psi a}$; and the correspondence of card-positions within the resulting hypercubes is now $\varphi:H\to H'$.  Thus we have produced a combinatorial theory that faithfully models the graphical notions of section 1.

Our next goal is to classify all hands up to 3 cards, including ``Sets''.

\section{All hand types with at most 3 cards}\label{secthreehands}

In figure \ref{foursets} we see representatives for four different types of hand.  We now extend this list to catalog \emph{every} type of hand with at most 3 cards.  We also compute the number of hands in each class.

\subsection*{Two-card hands}

\begin{figure}\centering\begin{tabular}{|lcr|}
	\hline
	Class & Common & Class\\
	representative & attributes & size\\
	\hline\\[-10pt]$\vcenter{\hbox{
	\includegraphics[width=.4in]{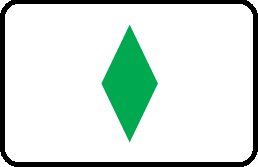}\,
	\includegraphics[width=.4in]{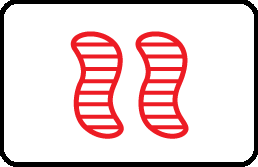}}}$ & None & 648
	\\[6pt]\hline\\[-10pt]$\vcenter{\hbox{
	\includegraphics[width=.4in]{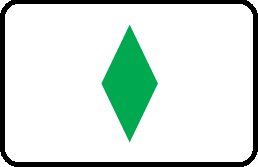}\,
	\includegraphics[width=.4in]{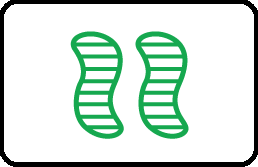}}}$ & One &  1,296
	\\[6pt]\hline\\[-10pt]$\vcenter{\hbox{
	\includegraphics[width=.4in]{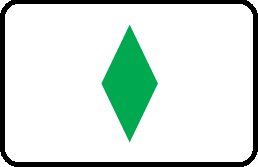}\,
	\includegraphics[width=.4in]{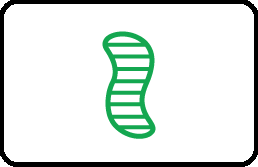}}}$ & Two &  972
	\\[6pt]\hline\\[-10pt]$\vcenter{\hbox{
	\includegraphics[width=.4in]{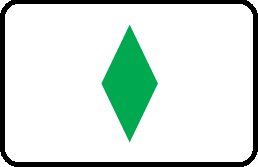}\,
	\includegraphics[width=.4in]{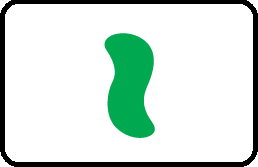}}}$ & Three &  324
	\\[6pt]\hline&&\\
	\textbf{Total:} & &3,240\\
	\hline
\end{tabular}\caption{All types of two-card hand.}\label{tabtwohands}\end{figure}

There are precisely four types of two-card hand.  Each type is distinguished by how many common attributes its cards agree in.  Figure \ref{tabtwohands} gives one example of each type, and counts how many hands are isomorphic to it.  For instance, there are 648 different two-card hands with no common attributes.

Computation of the numbers in figure \ref{tabtwohands} is straightforward.  For two-card hands with no common attributes, we first count 81 choices for the first card.  The second card must have either of two unused values in each attribute, for a total of $2^4=16$ choices.  This double-counts every hand (order doesn't matter) so we divide by two to get the answer:
\[81\cdot2^4/2=648.\]

Counting the two-card hands with one common attribute follows the same reasoning.  Now the second card must have an unused value in only three attributes, so the calculation becomes:
\[81\cdot{4\choose3}\cdot2^3/2=1296.\]

Similarly, for two common attributes, we calculate:
\[81\cdot{4\choose2}\cdot2^2/2=972.\]

The case of three common attributes is left as an exercise.

As a final check-sum, we observe that the total of all four class-sizes is 3,240, which equals $81\choose2$ as desired.

\subsection*{One-card hands}

\begin{figure}\centering\begin{tabular}{|lr|}
	\hline
	Class & Class\\
	representative & size\\
	\hline\\[-10pt]$\vcenter{\hbox{
	\includegraphics[width=.4in]{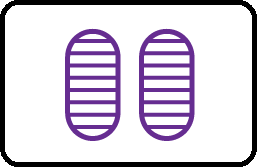}}}$ & 81
	\\[6pt]\hline&\\
	\textbf{Total:} & 81\\
	\hline
\end{tabular}\caption{There is only one type of one-card hand.}\label{tabonehand}\end{figure}

There is only one type of one-card hand.  All one-card hands are isomorphic.  See figure \ref{tabonehand}.  If $H=\{x\}$ and $H'=\{x'\}$ are one-card hands, then the map $\varphi (x)=x'$ is induced by taking (for instance) $\psi(a)=a$ and $\vartheta_a\left(v_a(x)\right)=v_a\left(x'\right)$ for each $a\in A$.

\subsection*{Zero-card hands}

The empty hand is unique, and thus represents its own class.  See figure \ref{tabzerohand}.

\begin{figure}\centering\begin{tabular}{|cc|}
	\hline
	Class & Class\\
	representative & size\\
	\hline
	$\emptyset$ & 1\\
	\hline&\\
	\textbf{Total:} & 1\\
	\hline
\end{tabular}\caption{There is only one type of zero-card hand.}\label{tabzerohand}\end{figure}

\subsection*{Three-card hands}

\begin{figure}\raisebox{12pt}{\begin{tabular}{|ccr|}
	\hline
	Class & & Class\\
	representative & Symbol & size\\
	\hline\\[-10pt]$\vcenter{\hbox{
	\includegraphics[height=1.5em]{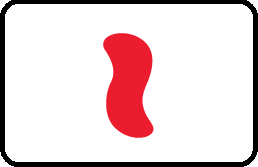}
	\includegraphics[height=1.5em]{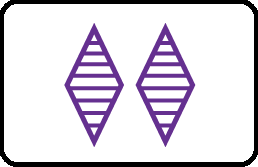}
	\includegraphics[height=1.5em]{81}}}$ &$(0;0,0,0)$& 216
	\\[3pt]\hline\\[-10pt]$\vcenter{\hbox{
	\includegraphics[height=1.5em]{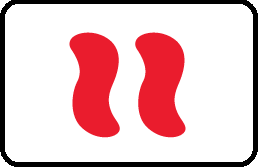}
	\includegraphics[height=1.5em]{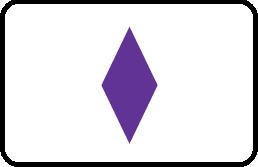}
	\includegraphics[height=1.5em]{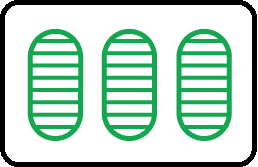}}}$ &$(0;0,0,1)$& 2,592
	\\[3pt]\hline\\[-10pt]$\vcenter{\hbox{
	\includegraphics[height=1.5em]{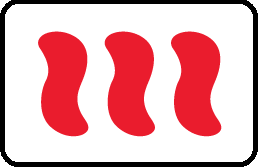}
	\includegraphics[height=1.5em]{4}
	\includegraphics[height=1.5em]{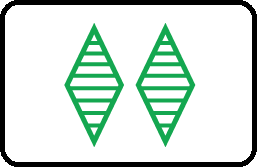}}}$ &$(0;0,0,2)$& 3,888
	\\[3pt]\hline\\[-10pt]$\vcenter{\hbox{
	\includegraphics[height=1.5em]{4}
	\includegraphics[height=1.5em]{1}
	\includegraphics[height=1.5em]{44}}}$ &$(0;0,0,3)$& 2,592
	\\[3pt]\hline\\[-10pt]$\vcenter{\hbox{
	\includegraphics[height=1.5em]{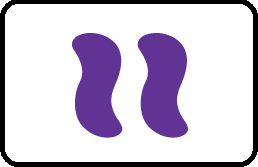}
	\includegraphics[height=1.5em]{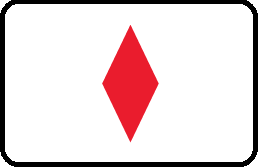}
	\includegraphics[height=1.5em]{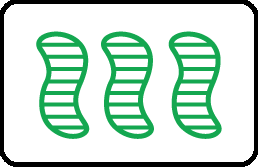}}}$ &$(0;0,1,1)$& 7,776
	\\[3pt]\hline\\[-10pt]$\vcenter{\hbox{
	\includegraphics[height=1.5em]{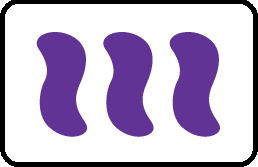}
	\includegraphics[height=1.5em]{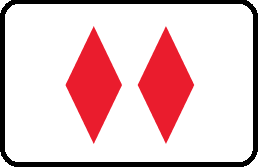}
	\includegraphics[height=1.5em]{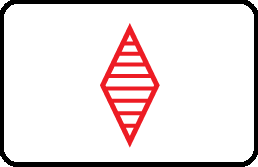}}}$ &$(0;0,1,2)$& 15,552
	\\[3pt]\hline\\[-10pt]$\vcenter{\hbox{
	\includegraphics[height=1.5em]{7}
	\includegraphics[height=1.5em]{11}
	\includegraphics[height=1.5em]{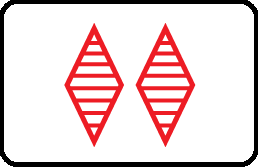}}}$ &$(0;0,1,3)$& 5,184
	\\[3pt]\hline\\[-10pt]$\vcenter{\hbox{
	\includegraphics[height=1.5em]{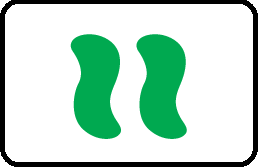}
	\includegraphics[height=1.5em]{1}
	\includegraphics[height=1.5em]{37}}}$ &$(0;0,2,2)$& 3,888
	\\[3pt]\hline\\[-10pt]$\vcenter{\hbox{
	\includegraphics[height=1.5em]{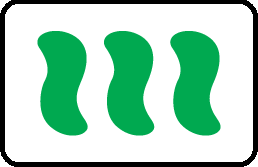}
	\includegraphics[height=1.5em]{10}
	\includegraphics[height=1.5em]{29}}}$ &$(0;1,1,1)$& 5,184
	\\[3pt]\hline\\[-10pt]$\vcenter{\hbox{
	\includegraphics[height=1.5em]{10}
	\includegraphics[height=1.5em]{9}
	\includegraphics[height=1.5em]{30}}}$ &$(0;1,1,2)$& 7,776
	\\[3pt]\hline
	\end{tabular}}
	\quad
	\begin{tabular}{|ccr|}
	\hline
	Class & & Class\\
	representative & Symbol & size\\
	\hline\\[-10pt]$\vcenter{\hbox{
	\includegraphics[height=1.5em]{11}
	\includegraphics[height=1.5em]{4}
	\includegraphics[height=1.5em]{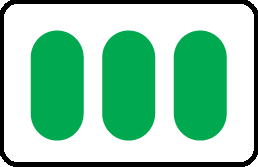}}}$ &$(1;0,0,0)$& 432
	\\[3pt]\hline\\[-10pt]$\vcenter{\hbox{
	\includegraphics[height=1.5em]{12}
	\includegraphics[height=1.5em]{1}
	\includegraphics[height=1.5em]{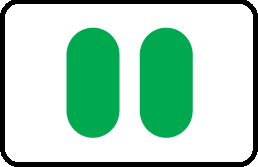}}}$ &$(1;0,0,1)$& 3,888
	\\[3pt] \hline\\[-10pt]$\vcenter{\hbox{
	\includegraphics[height=1.5em]{13}
	\includegraphics[height=1.5em]{4}
	\includegraphics[height=1.5em]{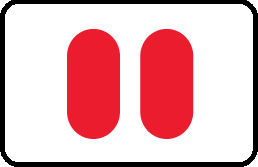}}}$ &$(1;0,0,2)$& 3,888
	\\[3pt]\hline\\[-10pt]$\vcenter{\hbox{
	\includegraphics[height=1.5em]{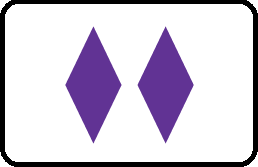}
	\includegraphics[height=1.5em]{1}
	\includegraphics[height=1.5em]{6}}}$ &$(1;0,1,1)$& 7,776
	\\[3pt]\hline\\[-10pt]$\vcenter{\hbox{
	\includegraphics[height=1.5em]{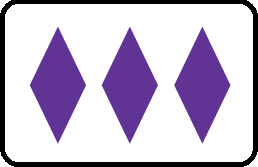}
	\includegraphics[height=1.5em]{1}
	\includegraphics[height=1.5em]{3}}}$ &$(1;0,1,2)$& 7,776
	\\[3pt]\hline\\[-10pt]$\vcenter{\hbox{
	\includegraphics[height=1.5em]{16}
	\includegraphics[height=1.5em]{1}
	\includegraphics[height=1.5em]{11}}}$ &$(1;1,1,1)$& 2,592
	\\[3pt]\hline\\[-10pt]$\vcenter{\hbox{
	\includegraphics[height=1.5em]{17}
	\includegraphics[height=1.5em]{2}
	\includegraphics[height=1.5em]{23}}}$ &$(2;0,0,0)$& 324
	\\[3pt]\hline\\[-10pt]$\vcenter{\hbox{
	\includegraphics[height=1.5em]{18}
	\includegraphics[height=1.5em]{3}
	\includegraphics[height=1.5em]{15}}}$ &$(2;0,0,1)$& 1,944
	\\[3pt]\hline\\[-10pt]$\vcenter{\hbox{
	\includegraphics[height=1.5em]{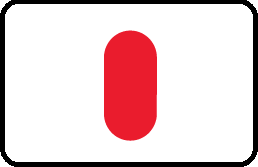}
	\includegraphics[height=1.5em]{1}
	\includegraphics[height=1.5em]{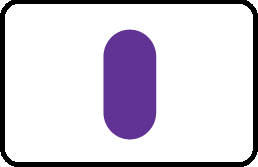}}}$ &$(2;0,1,1)$& 1,944
	\\[3pt]\hline\\[-10pt]$\vcenter{\hbox{
	\includegraphics[height=1.5em]{20}
	\includegraphics[height=1.5em]{23}
	\includegraphics[height=1.5em]{26}}}$ &$(3;0,0,0)$& 108
	\\[3pt]\hline&&\\
	\textbf{Total:} && 85,320\\
	\hline
\end{tabular}\caption{There are 20 types of three-card hand.}\label{tabthreehands}\end{figure}

There are 20 types of three-card hand, shown in figure \ref{tabthreehands}.  They are classified according to their \emph{symbols}, as follows.

\begin{definition}
If $H$ is a three-card hand, its \emph{symbol} is a 4-tuple of numbers $(t;p_1,p_2,p_3)$, where the $p_i$'s are unordered.  Here, $t$ is the number of attributes common to all three cards, and the $p_i$'s count the common attributes of each pair, besides those already counted by $t$.
\end{definition}

\begin{example}
The hand \includegraphics[height=1.5em]{45}\,\includegraphics[height=1.5em]{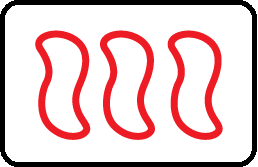}\,\includegraphics[height=1.5em]{66} has symbol $(1;0,1,2)$ because there is one attribute (number) common to all three cards, and there are 0, 1, \& 2 other attributes common to pairs of cards.
\end{example}

\begin{example}
A three-card hand is a ``Set'' if its symbol is $(t;0,0,0)$ with $t\in\{0,1,2,3\}$.  The ``Sets'' of figure \ref{foursets}, for instance, have symbols $(0;0,0,0)$, $(1;0,0,0)$, $(2;0,0,0)$, and $(3;0,0,0)$ (respectively as listed from top to bottom).
\end{example}

The next two theorems show that the classification in figure \ref{tabthreehands} is accurate and complete.

\begin{theorem}
Two three-card hands are isomorphic if and only if they have the same symbol.

\proof $\Longleftarrow$:  Suppose $H$ and $H'$ both have symbol $(t;p_1,p_2,p_3)$.  The numbers in this symbol count certain attributes shared by cards in $H$; they also count attributes shared by cards in $H'$.  Let $\psi:A\to A$ be the permutation sending the $t$ common attributes of $H$ to the those of $H'$; sending the $p_1$ attributes shared by a pair of cards in $H$ to the $p_1$ attributes shared by a pair of cards in $H'$; and likewise for $p_2$ and $p_3$.  Value bijections $\vartheta_a$ can now be chosen for each attribute so that the $a$-values present in $H$ map to corresponding $\psi a$-values in $H'$.  This induces an isomorphism $\varphi:H\to H'$.

$\Longrightarrow$:  Suppose $\varphi:H\to H'$ is an isomorphism induced by some permutation $\psi:A\to A$ and bijections $\vartheta_a:V_a\to V_{\psi a}$ for each $a\in A$.  Suppose $H$ has symbol $(t;p_1,p_2,p_3)$ and $H'$ has symbol $(t';p_1',p_2',p_3')$.  If $a$ is an attribute common to all three cards of $H$, then $\psi a$ is an attribute common to all three cards of $H'$.  Since $\psi$ is one-to-one, there must be the same number of each, that is, $t=t'$. Now consider only the attributes not yet counted. If attribute $b$ is common to a particular pair of cards in $H$, then attribute $\psi b$ is common to the corresponding pair of cards in $H'$.  If that pair in $H$ has $p_1$ attributes in common, then the corresponding pair in $H'$ also has $p_1$ attributes in common, so $p_1$ is among the $p_i'$'s.  The same applies to $p_2$ and $p_3$.  Thus $H$ and $H'$ have the same symbol.  \quad$\square$
\end{theorem}

\begin{theorem}
A 4-tuple $(t;p_1,p_2,p_3)$ is the symbol of some three-card hand if and only if 
\begin{align*}
	t				&\in\{0,1,2,3\}\\
	\text{each }p_i	&\in\{0,1,2,3\}\\
	\text{each }p_i+t	&\in\{0,1,2,3\}\\
	t+p_1+p_2+p_3	&\in\{0,1,2,3,4\}.
\end{align*}

\proof  $\Longrightarrow$:  Each number $t,p_1,p_2,p_3$ counts attributes shared by two or more cards.  Two cards can agree in at most three attributes, so none of these numbers can exceed 3.  This gives us the first two conditions.  The sum $p_i+t$ counts the \emph{total} number of attributes shared by two cards; again, this is at most three.  That gives us the third condition.  Lastly, the numbers $t,p_1,p_2,p_3$ each counts \emph{different} attributes--- each attribute is counted at most once.  There are only four attributes total, so the sum $t+p_1+p_2+p_3$ is at most four.  This gives us the last condition.  Thus, each of the listed conditions is necessary.

$\Longleftarrow$:  In figure \ref{tabthreehands}, every possible tuple satisfying these conditions is systematically listed.  The table also gives a three-card hand having each such tuple as its symbol.  Thus, the four conditions are sufficient. \quad$\square$
\end{theorem}

The class sizes in figure \ref{tabthreehands} were computed as follows.  For each symbol $(t;p_1,p_2,p_3)$ we construct a hand $H=\{x_1,x_2,x_3\}$ that has $t$ attributes shared all-around, plus $p_1$ shared by $x_2$ and $x_3$, $p_2$ shared by $x_1$ and $x_3$, and $p_3$ shared by $x_1$ and $x_2$.  Then:
\begin{equation}\label{threehandformula}\vcenter{\hbox{
\includegraphics[width=4in]{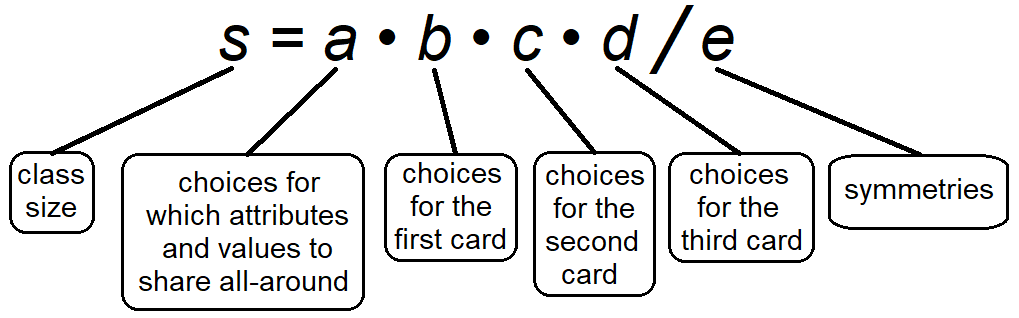}\\
}}\end{equation}
These choices are made in order, from left to right.  Each factor is determined by the symbol $(t;p_1,p_2,p_3)$.  The first two factors, $a$ and $b$, depend only on $t$:
\begin{itemize}
	\item If $t=0$, then $a=1$ and $b=81$.
	\item If $t=1$, then $a=12$ and $b=27$.
	\item If $t=2$, then $a=54$ and $b=9$.
	\item If $t=3$, then $a=108$ and $b=3$.
\end{itemize}
The next two factors, $c$ and $d$, depend more subtly on the symbol, as illustrated in the examples below.  The number $e$ counts \emph{symmetries} of the hand $H$, that is, the self-isomorphisms of $H$.  This depends on how many distinct values are among the $p_i$:
\begin{itemize}
	\item If $p_1=p_2=p_3$, then $e=6$.
	\item If there are exactly two distinct values among the $p_i$, then $e=2$.
	\item If the $p_i$ are all different, then $e=1$.
\end{itemize}

\begin{example}
How many three-card hands have symbol $(0;0,0,0)$?  There is only $a=1$ way to choose zero attributes to be shared all-around.  There are $b=81$ possible choices for the first card.  With the first card fixed, there are $c=16$ cards that have nothing in common with it; choose one of these as the second card.  There is only $d=1$ card that has nothing in common with either of the first two chosen cards.  In our symbol, $p_1=p_2=p_3$, so there are $e=6$ symmetries for this hand.  Plugging these values into equation (\ref{threehandformula}) yields the answer:  There are $1\cdot81\cdot16\cdot1/6 = 216$ hands with this symbol.
\end{example}

\begin{example}
How many three-card hands have symbol $(1;0,1,2)$?  One attribute will be shared all-around; there are 4 choices of which attribute and 3 choices for its value, so $a=12$.  With that choice fixed, there are $b=27$ choices for the first card.  The second card should have two additional things in common with the first; there are $c=6$ such cards.  The third card must match the first in one additional way, and must have nothing else in common with the second; there are $d=4$ such cards.  There are no repeated values among the $p_i$ in our symbol, so $e=1$.  By equation (\ref{threehandformula}), there are $12\cdot27\cdot6\cdot4/1=7,776$ hands with this symbol.
\end{example}

Following these examples, the reader is invited to try computing the other class sizes in figure \ref{tabthreehands} as an exercise.  In particular, verify that the four ``Sets'' in figure \ref{tabthreehands} have class sizes 216, 432, 324, and 108.  The sum of these numbers is 1080, the total number of ``Sets'' that exist, as observed in section \ref{intro}.  The sum of \emph{all} the class sizes in figure \ref{tabthreehands} is 85,320, which equals $81\choose3$ as expected.

\section{The general hand-classification problem}

In section \ref{secthreehands} we classified small subsets of the standard 81-card \set\, deck.  This gives the first few entries in figure \ref{tabcountclasses}. Up to isomorphism, there are four types of two-card hand (figure \ref{tabtwohands}), 20 types of three-card hand (figure \ref{tabthreehands}), and so forth. The first table entry not computed in section \ref{secthreehands} states that there are 144 types of four-card hand; our present goal is to complete the entire table.

The numbers in figure \ref{tabcountclasses} appear in the OEIS as integer sequence \href{http://www.oeis.org/A034216}{\underline{A034216}}, the ``number of ternary codes of length 4 with $n$ words'' \cite{OEIS}.  Evidently, every \set\, card can be thought of as a 4-letter word in a ternary code, and a \set\, hand with $n$ cards can be thought of a code with $n$ words.  Then, \set\, hands are isomorphic exactly when they correspond to isometric ternary codes.\footnote{Thank you to Andrew Sward and Jordan Thompson for this insight.}

The notion of ternary code isometry is perhaps a bit abstruse for this study. We would like to verify all the entries in figure \ref{tabcountclasses} using the methods shown in section \ref{secthreehands}.  To do this, we need to devise a ``symbol'' for $n$-card hands, as we did for $n=3$ in figure \ref{tabthreehands}.  This is left open as a challenge to the reader. How, for example, can you quickly recognize whether two four-card hands are isomorphic?

\begin{figure}\centering\begin{tabular}{r|c|c|c|c|c| c |c|c|c|c|c|}
Hand size &0&1&2&3&4&$\cdots$&77&78&79&80&81\\
\hline
\# classes &1&1&4&20&144&$\cdots$&144&20&4&1&1
\end{tabular}\caption{How many types of $n$-card hand are there?}\label{tabcountclasses}\end{figure}

The next theorem says that figure \ref{tabcountclasses} is a palindrome, so there are only 41 cases to work out, rather than 82.
\begin{theorem}  Two hands $H$ and $H'$ are isomorphic if and only if their complements are isomorphic.  Therefore, the number of classes of $n$-card hands equals the number of classes of $(81-n)$-card hands.

\proof Suppose $\varphi:H\to H'$ is an isomorphism induced by the bijections $\psi$ and $\vartheta_a$ for all $a\in A$.  Then $\varphi$ extends to a self-isomorphism the entire deck (induced by the same $\psi$ and $\vartheta_a$).  This self-isomorphism restricts to an isomorphism of the complements $H^\mathcal C\to H'^{\mathcal C}$.\quad$\square$
\end{theorem}

Thus far, we have worked with the standard \set\, deck, having 4 attributes with 3 values each.  The ideas we have studied generalize to other decks, where these numbers may be different.

\begin{definition}
A \emph{\set-style deck}, denoted $D(k^d)$, is a collection of $k^d$ cards, each card distinguished by one of $k$ values in each of $d$ attributes.  The standard 81-card \set\, deck is $D(3^4)$.
\end{definition}

The \textbf{general hand-classification problem} asks three questions.
\begin{enumerate}
	\item {\bfseries In the deck $D(k^d)$, given $n\in\mathbb N$, how many types of $n$-card hand exist?}\\
	Figure \ref{tabcountclasses} answers this question for small hands in a standard \set\, deck. But how can we produce a similar table for larger decks?

	 The question is equivalent to:  How many $k$-nary codes of length $d$ with $n$ words exist? 

	As already noted, length-4 ternary codes correspond to ordinary \set\, cards; this correspondence generalizes to larger decks. For some values of $k$, $d$, and $n$, the answer can be found in the OEIS. For example, according to sequence \href{http://www.oeis.org/A034240}{\underline{A034240}} in the OEIS\footnote{``Number of quaternary codes of length 9 with n words.'' \cite{OEIS2}}, there are approximately $1.08\times10^{34}$ classes of eleven-card hands from the \set-style deck with nine attributes and four values per attribute, a.k.a. $D(4^9)$. Other sequences in the OEIS correspond to other decks, too.

	A detailed analysis of $k$-nary codes and their relation to general \set\, decks will be a focus of future work on this project.
	\item {\bfseries How large is the isomorphism class of a given a hand $H$ in $D(k^d)$?}\\
	 This question asks for the likelihood of drawing a particular type of hand from a shuffled deck-- useful information when designing a game.
	\item {\bfseries How can one quickly determine whether two given hands $H,H'\subset D(k^d)$ are isomorphic?}\\
	 This question requests an efficient general algorithm for comparing two hands.  For three-card hands of the standard \set\, deck, the algorithm is simply the ``symbol'' used in figure \ref{tabthreehands}.  But a general method for larger hands and larger decks is not known.
\end{enumerate}

Equation (\ref{threehandformula}) makes reference to the \emph{symmetries} of a three-card hand, that is, the self-isomorphisms of the hand.  In general, the symmetries of an $n$-card hand $H$ form a subgroup $\Aut(H)$ of $S_H$, the group of all $n!$ permutations of $H$.  The group $\Aut(H)$ depends only on the isomorphism type of $H$, so it may be used as a distinguishing type invariant for hands.

\begin{example}\label{symmetricfour}
From the standard \set\, deck, let $H=\includegraphics[height=1.5em]{10}\,\includegraphics[height=1.5em]{20}\,\includegraphics[height=1.5em]{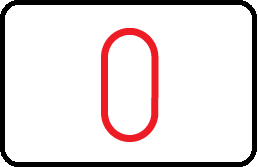}\,\includegraphics[height=1.5em]{22}$\,.  All 24 permutations of this hand are symmetries.  To see this, note that each card differs from all the others in one attribute.  Suppose $\varphi:H\to H$ is a permutation.  For each card $x\in H$, let $\psi$ map the attribute $a$ distinguishing $x$ to the attribute $\psi a$ distinguishing $\varphi x$.  Let $\vartheta_a$ map the value $v_a(x)$ to $v_{\psi a}(\varphi x)$.  These assignments produce bijections that induce $\varphi$ as an isomorphism.  Thus, $\Aut(H)=S_H$.
\end{example}

\begin{example}
By contrast, the hand $H=\includegraphics[height=1.5em]{10}\,\includegraphics[height=1.5em]{20}\,\includegraphics[height=1.5em]{73}\,\includegraphics[height=1.5em]{19}$ has only six symmetries.  To see this, suppose $\varphi:H\to H$ is a permutation fixing \includegraphics[height=1em]{19}\,.  Let $\psi$ and $\vartheta_a$ be as in example \ref{symmetricfour}; note that now, $\psi(\text{\small{color}})=\text{\small{color}}$.  As in example \ref{symmetricfour}, these bijections induce $\varphi$ as an isomorphism.  On the other hand, suppose $\varphi':H\to H$ is a permutation \emph{not} fixing \includegraphics[height=1em]{19}\,.  Since \includegraphics[height=1em]{19}\, has three attributes in common with each other card, while every other pair of cards has only two common attributes, $\varphi'$ cannot be an isomorphism.
\end{example}

The symmetries of a given hand are generally not the same as the symmetries of its complement, as the following examples demonstrate.

\begin{example}
The entire \set\, deck has 31,104 symmetries \shortcite{McMahon}.  Every choice of bijections $\psi,\vartheta_a$ induces a unique isomorphism of the deck.  There are $4!$ choices for $\psi$, and $3!$ choices for each of the four maps $\vartheta_a$.  Thus
\[
|\Aut(\text{\small{the deck}})|=4!\cdot3!^4=31,104.
\]
The empty hand, of course, has only one symmetry, which is equally induced by any of these choices of bijections.
\end{example}

\begin{example}
Let $H$ be the entire \set\, deck except for the single card \includegraphics[height=1.5em]{26}\,.   Regardless of $\psi$, for each attribute $a$, the bijection $\vartheta_a$ must take $v_a( \vcenter{\hbox{\includegraphics[height=1em]{26}}})$ to $v_{\psi a}(\vcenter{\hbox{\includegraphics[height=1em]{26}}})$.  Thus there are still $4!$ choices for $\psi$, but there are now only 2 choices for each $\vartheta_a$.  Thus
\[
|\Aut(H)|=4!\cdot2^4=384.
\]
The hand $H^\mathcal C=\includegraphics[height=1.5em]{26}$ has only one symmetry, which is equally induced by any choice of bijections satisfying the above condition.
\end{example}

The preceding examples also demonstrate how several different choices of $\psi$ and $\vartheta_a$ may induce the same isomorphism of $H$.  Finding all possible choices inducing a given $\varphi$ is another good combinatorial problem:

\begin{enumerate}
\item[4.] {\bfseries Given an isomorphism $\varphi$ between hands, list all possible choices of attribute- and value-correspondences inducing $\varphi$.}\\
The table in example \ref{exiso} shows one possible set of correspondences that induce the given isomorphism $\varphi$. Is this table forced, or are there other possibilities?
\end{enumerate}

\section{\stun\, and other games}\label{gamesection}

We now present a variant of the game \set, which we call \stun\footnote{The name was suggested by a friend as a near-anagram of \emph{UN-SET}.}.

First, we review the rules of \set.  From the standard \set\, deck, deal twelve cards face-up on the table.  All players now search for a ``Set'', that is, three cards having either one or three values--- but not precisely two--- in each attribute.  If no ``Set'' is present, more cards are dealt.  When the deck is exhausted, the player who collected the most ``Sets'' wins.

The rules of \stun\, are the same except in two details:  Now, only nine cards initially should be dealt (rather than twelve), and instead of collecting ``Sets'', players now collect ``Stuns''.

\begin{definition}  A ``Stun'' is three cards having precisely \emph{two} values in each attribute.  \end{definition}

\begin{figure}[H]\centering
\includegraphics[width=0.5in]{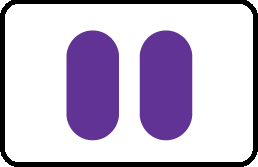}
\includegraphics[width=0.5in]{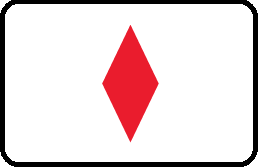}
\includegraphics[width=0.5in]{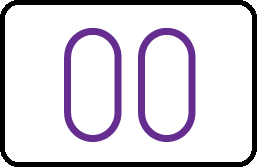}
\caption{A ``Stun''.}\label{astun}
\end{figure}

\begin{figure}\centering\begin{tabular}{ccc}
\includegraphics[width=0.5in]{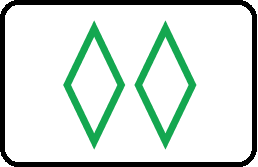}&
\includegraphics[width=0.5in]{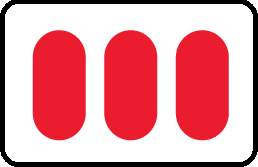}&
\includegraphics[width=0.5in]{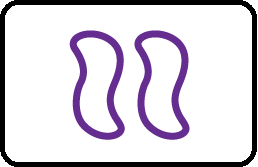}\\
\includegraphics[width=0.5in]{38}&
\includegraphics[width=0.5in]{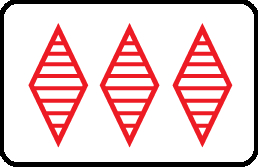}&
\includegraphics[width=0.5in]{18}\\
\includegraphics[width=0.5in]{8}&
\includegraphics[width=0.5in]{30}&
\includegraphics[width=0.5in]{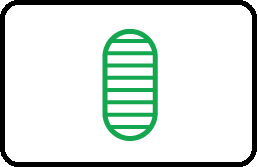}
\end{tabular}\caption{There are twenty ``Stuns'' here.  How many can you find?}\label{stungame}\end{figure}

See figure \ref{stungame} for a sample gameboard.  Of the ${9\choose3}=84$ possible three-card hands, fully 20 of them are ``Stuns''.  How many can you find?

\begin{figure}\centering\begin{tabular}{|lcr|}
	\hline
	Class		&			& Class\\
	representative & Symbol	& size\\
	[3pt]\hline\\[-10pt]$\vcenter{\hbox{
	\includegraphics[height=1.5em]{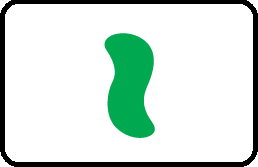}
	\includegraphics[height=1.5em]{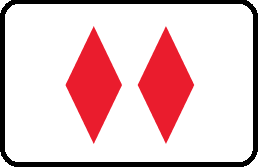}
	\includegraphics[height=1.5em]{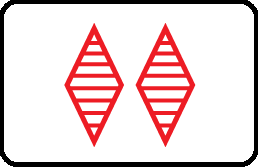}}}$ &$(0;0,1,3)$& 5,184
	\\[3pt]\hline\\[-10pt]$\vcenter{\hbox{
	\includegraphics[height=1.5em]{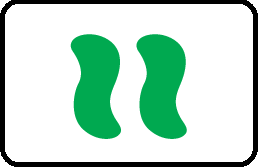}
	\includegraphics[height=1.5em]{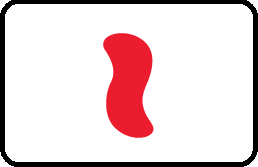}
	\includegraphics[height=1.5em]{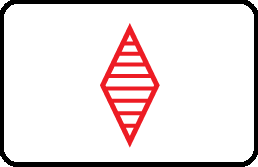}}}$ &$(0;0,2,2)$& 3,888
	\\[3pt]\hline\\[-10pt]$\vcenter{\hbox{
	\includegraphics[height=1.5em]{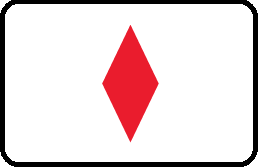}
	\includegraphics[height=1.5em]{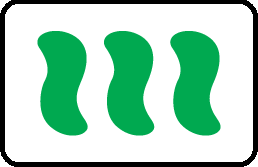}
	\includegraphics[height=1.5em]{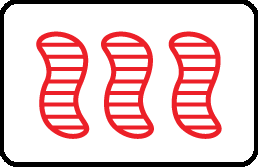}}}$ &$(0;1,1,2)$& 7,776
	\\[6pt]\hline&&\\
	\textbf{Total:} & &16,848\\
	\hline
\end{tabular}\caption{The three types of ``Stun''.}\label{tabstuns}\end{figure}

Whereas ``Sets'' come in four varieties (figure \ref{foursets}), ``Stuns'' come in three.  The symbol $(t;p_1,p_2,p_3)$ for a ``Stun'' must have $t=0$, since no value is shared all-around, and must satisfy $p_1+p_2+p_3=4$, since each of the four attributes is common to one pair of cards.  Only three entries in figure \ref{tabthreehands} satisfy these requirements; these are distilled in figure \ref{tabstuns}.

There are 16,848 ``Stuns'' in the deck, out of 85,320 possible three-card hands.  The probability that the top three cards of a shuffled deck are a ``Stun'' is $16848/85320\approx0.1975$--- a nearly one-in-five chance!  Contrast this with the one-in-79 chance of drawing a ``Set'', and it would seem that, by the numbers, \stun\, is an inferior game to \set.

Nonetheless, we have heartily enjoyed playing \stun.  The challenge of instantly recognizing a ``Stun'' is an exciting new skill to learn.  It's especially fun to play \stun\, against seasoned \set-players, who struggle to suppress their conditioned ``Set''-grabbing reflexes.  For an extra laugh, play alternating games of \set\, and \stun\, without shuffling the deck in between, and watch your opponents try to switch their brains back-and-forth between the two modes.

For \stun\, bonus points, try to have no cards left over at the end of the game.  It's an amusing puzzle to partition the last nine cards into three ``Stuns'', when possible (which it usually is).

Many interesting combinatorial questions about \set\, can also be asked about \stun.  For instance, it's well-known that 21 cards suffice to guarantee the presence of a ``Set''.  How many cards are needed to guarantee a ``Stun''?  The answer is at least 28 (the 27 purple cards, for instance, contain no ``Stun'').

The idea behind \set\, and \stun\, can be extrapolated to devise other games in a similar vein.  Certain isomorphism classes of hand are designated as goals, and players search among some face-up cards for specimens.  Three-card-goal games may be invented by choosing one or more classes from figure \ref{tabthreehands}.  The game of \emph{SOOT}, for example, might require players to collect hands of type (0;0,1,2), such as 
	\includegraphics[height=1em]{6}
	\includegraphics[height=1em]{11}
	\includegraphics[height=1em]{37}\,.  While it may be very challenging to recognize hands of this type, there are 15,552 of them in the deck, so the probability of finding one by accident is rather high ($\approx0.1823$).  The ubiquity of the goal offsets the perceptual challenge of recognizing it, yielding a balanced game.  Exercise:  How many ``Soots'' can you find in figure \ref{stungame}?
	
\begin{figure}[t]\centering
\begin{tabular}{cc|ccccc}
\includegraphics[height=2em]{25}&\hspace{\fill}&\hspace{\fill}&
\includegraphics[height=2em]{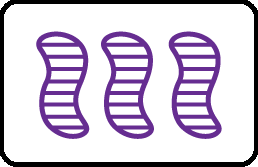}&
\includegraphics[height=2em]{12}&
\includegraphics[height=2em]{54}&
\includegraphics[height=2em]{2}\\

\includegraphics[height=2em]{79}&&&
\includegraphics[height=2em]{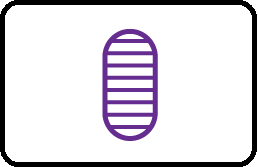}&
\includegraphics[height=2em]{23}&
\includegraphics[height=2em]{26}&
\includegraphics[height=2em]{1}\\

\includegraphics[height=2em]{3}&&&
\includegraphics[height=2em]{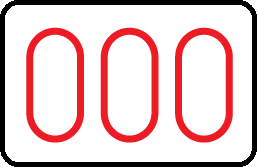}&
\includegraphics[height=2em]{66}&
\includegraphics[height=2em]{21}&
\includegraphics[height=2em]{59}\\

\includegraphics[height=2em]{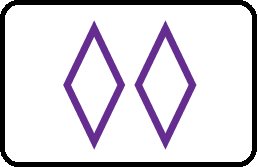}&&&
\includegraphics[height=2em]{20}&
\includegraphics[height=2em]{9}&
\includegraphics[height=2em]{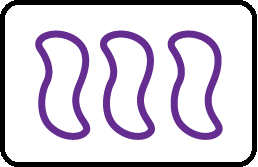}&
\includegraphics[height=2em]{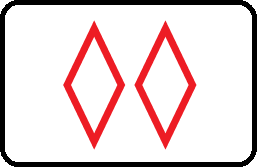}\\

\Large{Goal} && \multicolumn{4}{c}{\Large{Playing area}}
\end{tabular}
\caption{Find four cards in the playing area isomorphic to the goal.  (Solution in figure \ref{solution}.)}\label{hardgame}\end{figure}

For the truly ambitious player, larger hands may also be used as goals.  Shuffle the deck and deal four cards off to the side--- this four-card hand type is your goal.  Now deal sixteen cards into the playing area and search for four among them that are isomorphic to the goal.  See figure \ref{hardgame}.  Mathematical analysis of this game would be a computational challenge, but would follow the essential principles laid out in this paper.

\begin{figure}
	\centering
	\includegraphics[width=3in]{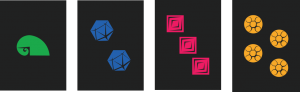}
	\caption{A ``quad'' in the game EvenQuads.\\
		 Image borrowed from \protect\cite{AWM}.}
	\label{evenquads}
\end{figure}

Recently, the AWM released a novel \set-style deck called \emph{EvenQuads} \cite{AWM}, The cards have four values in each of three attributes. That is, the EvenQuads deck is $D(4^3)$ with 64 cards. The game is played similarly to \set, though now the objective is to find "quads", which are four cards that, in each attribute, are either all alike, all different, or are split 2-2. (An example appears in figure \ref{evenquads}.) By definition, the set of "Quads" is a union of isomorphism classes of four-card hands. Thus EvenQuads lends itself to the same sort of analysis and variations as laid out in this paper.

\pagebreak
\section{Further reading}

Many excellent treatments have been written about the mathematics of \set, including a few that define equivalence relations on hands \cite{Coleman, McMahon}.  However, the equivalence we use here has not been described previously, that I can see.  There is a reference to ``isomorphisms'' of hands in the documentation for a 2001 CWEB program by Donald Knuth \cite{Knuth}. This program sought to classify, up to isomorphism, maximal hands that include no ``Sets''.  The documentation may be read \href{https://www-cs-faculty.stanford.edu/~knuth/programs.html}{\underline{here}}. (Scroll down to find ``SETSET''.)\\

\begin{figure}[p]\centering
\begin{tabular}{cc|ccccc}
\includegraphics[height=2em]{25}&\hspace{\fill}&\hspace{\fill}&
\includegraphics[height=2em]{33}&
\includegraphics[height=2em]{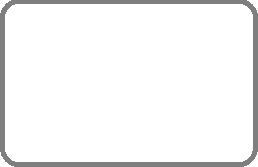}&
\includegraphics[height=2em]{blank}&
\includegraphics[height=2em]{blank}\\

\includegraphics[height=2em]{79}&&&
\includegraphics[height=2em]{blank}&
\includegraphics[height=2em]{blank}&
\includegraphics[height=2em]{26}&
\includegraphics[height=2em]{blank}\\

\includegraphics[height=2em]{3}&&&
\includegraphics[height=2em]{blank}&
\includegraphics[height=2em]{66}&
\includegraphics[height=2em]{blank}&
\includegraphics[height=2em]{blank}\\

\includegraphics[height=2em]{68}&&&
\includegraphics[height=2em]{blank}&
\includegraphics[height=2em]{blank}&
\includegraphics[height=2em]{blank}&
\includegraphics[height=2em]{65}\\

\Large{Goal} &&& \multicolumn{4}{c}{\Large{Playing area}}\\
\hline
\end{tabular}
\[
\begin{matrix}
\text{Goal}=	&\vcenter{\hbox{\includegraphics[width=.4in]{25}}}
	&\vcenter{\hbox{\includegraphics[width=.4in]{79}}}
	&\vcenter{\hbox{\includegraphics[width=.4in]{3}}}
	&\vcenter{\hbox{\includegraphics[width=.4in]{68}}}\\
\varphi:&\mapsdown&\mapsdown&\mapsdown&\mapsdown\\
\text{Answer}=	&\vcenter{\hbox{\includegraphics[width=.4in]{65}}}
	&\vcenter{\hbox{\includegraphics[width=.4in]{66}}}
	&\vcenter{\hbox{\includegraphics[width=.4in]{26}}}
	&\vcenter{\hbox{\includegraphics[width=.4in]{33}}}
\end{matrix}
\]
\begin{tabular}{|rcl|rcl|}
	\hline
	\textbf{COLOR}	&$\Longmapsto$& \textbf{COLOR}	&\textbf{FILL}	&$\Longmapsto$& \textbf{NUMBER}\\
	red				&$\mapsto$& 	green 			& empty			&$\mapsto$& 	triple\\
	green			&$\mapsto$& 	red				& stripe			&$\mapsto$& 	single\\
	purple			&$\mapsto$& 	purple			& solid			&$\mapsto$& 	double\\
	\hline
	\textbf{NUMBER}	&$\Longmapsto$&\textbf{SHAPE}	&\textbf{SHAPE}	&$\Longmapsto$&\textbf{FILL}\\
	single			&$\mapsto$& 	diamond 		& diamond		&$\mapsto$& 	stripe\\
	double			&$\mapsto$& 	squiggle			& oval			&$\mapsto$& 	empty\\
	triple			&$\mapsto$& 	oval				& squiggle		&$\mapsto$& 	solid\\
	\hline	
\end{tabular}

\caption{One possible solution to figure \ref{hardgame}.}\label{solution}\end{figure}

\bibliographystyle{apacite}
\raggedright\bibliography{SET}
\end{document}